\documentclass[12pt]{amsart}
\def\Image{\operatorname{Im}}
\def\Min{\operatorname{Min}}
\def\N{{\Bbb N}}

\def\R{{\Bbb R}}
\newtheorem{Theorem}{Theorem}[section]
\newtheorem{Corollary}[Theorem]{Corollary}

\theoremstyle{definition}

\newtheorem{Problem}[Theorem]{Problem}

\theoremstyle{remark}

\begin{document}
\sloppy
\title{On fixed-point sets in the boundary of a CAT(0) space}
\author{Tetsuya Hosaka} 
\address{Department of Mathematics, Utsunomiya University, 
Utsunomiya, 321-8505, Japan}
\date{May 2, 2005}
\email{hosaka@cc.utsunomiya-u.ac.jp}
\thanks{
Partly supported by the Grant-in-Aid for Young Scientists (B), 
The Ministry of Education, Culture, Sports, Science and Technology, Japan.
(No.~15740029).}
\maketitle
\begin{abstract}
In this paper, 
we investigate the fixed-point set of an element of a CAT(0) group in its boundary.
Suppose that a group $G$ acts geometrically on a CAT(0) space $X$.
Let $g\in G$ and 
let $\mathcal{F}_g$ be the fixed-point set of $g$ in the boundary $\partial X$.
Then we show that $\mathcal{F}_g=L(Z_g)$, 
where $Z_g$ is the centralizer of $g$ 
(i.e.\ $Z_g=\{v\in G|\, gv=vg\}$) and 
$L(Z_g)$ is the limit set of $Z_g$ in $\partial X$.
Thus we obtain that 
$\mathcal{F}_g\neq \emptyset$ if and only if 
the set $Z_g$ is infinite.
We also show that if $g$ is a hyperbolic isometry, 
then $\mathcal{F}_g=\partial\Min(g)$, 
where $\partial\Min(g)$ is 
the boundary of the minimal set $\Min(g)$ of $g$.
This implies that the fixed-point set $\mathcal{F}_g$ and 
the periodic-point set $\mathcal{P}_g$ of $g$ in $\partial X$ 
have suspension forms.
\end{abstract}

%%%%%%%%%%%%%
% Section 1 %
%%%%%%%%%%%%%
\section{Introduction and preliminaries}

The purpose of this paper is to study 
the fixed-point set of an element of a CAT(0) group in its boundary.

We say that a metric space $(X,d)$ is a {\it geodesic space} if 
for each $x,y \in X$, 
there exists an isometric embedding $\xi:[0,d(x,y)] \rightarrow X$ such that 
$\xi(0)=x$ and $\xi(d(x,y))=y$ (such $\xi$ is called a {\it geodesic}).
Also a metric space $(X,d)$ is said to be {\it proper} 
if every closed metric ball is compact.

Let $(X,d)$ be a geodesic space and 
let $T$ be a geodesic triangle in $X$.
A {\it comparison triangle} for $T$ is 
a geodesic triangle $\overline{T}$ in the Euclidean plane $\R^2$
with same edge lengths as $T$.
Choose two points $x$ and $y$ in $T$. 
Let $\bar{x}$ and $\bar{y}$ denote 
the corresponding points in $\overline{T}$.
Then the inequality $$d(x,y) \le d_{\R^2}(\bar{x},\bar{y})$$ 
is called the {\it CAT(0)-inequality}, 
where $d_{\R^2}$ is the natural metric on $\R^2$.
A geodesic space $(X,d)$ is called a {\it CAT(0) space} 
if the CAT(0)-inequality holds
for all geodesic triangles $T$ and for all choices of two points $x$ and 
$y$ in $T$.

Let $X$ be a proper CAT(0) space and $x_0 \in X$.
The {\it boundary of $X$ with respect to $x_0$}, 
denoted by $\partial_{x_0}X$, is defined as 
the set of all geodesic rays issuing from $x_0$. 
Then we define a topology on $X \cup \partial_{x_0}X$ 
by the following conditions: 
\begin{enumerate}
\item[(1)] $X$ is an open subspace of $X \cup \partial_{x_0}X$. 
\item[(2)] For $\alpha \in \partial_{x_0}X$ and $r, \epsilon >0$, let
$$ U_{x_0}(\alpha;r,\epsilon)=
\{ x \in X \cup \partial_{x_0} X \,|\, 
x \not\in B(x_0,r),\ d(\alpha(r),\xi_{x}(r))<\epsilon \}, $$
where $\xi_{x}:[0,d(x_0,x)]\rightarrow X$ is the geodesic from $x_0$ to $x$
($\xi_{x}=x$ if $x \in \partial_{x_0} X$).
Then for each $\epsilon_0>0$, 
the set 
$$\{U_{x_0}(\alpha;r,\epsilon_0)\,|\, r>0\}$$ 
is a neighborhood basis for $\alpha$ in $X \cup \partial_{x_0}X$. 
\end{enumerate}
This is called the {\it cone topology} on $X \cup \partial_{x_0} X$.
It is known that 
$X \cup \partial_{x_0} X$ is 
a metrizable compactification of $X$ (\cite{BH}, \cite{GH}).

Let $X$ be a geodesic space.
Two geodesic rays $\xi, \zeta:[0,\infty) \rightarrow X$ are 
said to be {\it asymptotic} if there exists a constant $N$ such that 
$d(\xi(t),\zeta(t)) \le N$ for each $t \ge 0$. 
It is known that 
for each geodesic ray $\xi$ in $X$ and 
each point $x \in X$, 
there exists a unique geodesic ray $\xi'$ issuing from $x$ 
such that $\xi$ and $\xi'$ are asymptotic.

Let $x_0$ and $x_1$ be two points of a proper CAT(0) space $X$.
Then there exists a unique bijection 
$\Phi:\partial_{x_0}X \rightarrow \partial_{x_1}X$ 
such that $\xi$ and $\Phi(\xi)$ are asymptotic 
for each $\xi \in \partial_{x_0}X$. 
It is known that $\Phi:\partial_{x_0}X\rightarrow \partial_{x_1}X$ 
is a homeomorphism (\cite{BH}, \cite{GH}).

Let $X$ be a proper CAT(0) space.
The asymptotic relation is an equivalence relation 
in the set of all geodesic rays in $X$.
The {\it boundary of} $X$, denoted by $\partial X$, 
is defined as the set of asymptotic equivalence classes of geodesic rays.
The equivalence class of a geodesic ray $\xi$ is denoted by $\xi(\infty)$.
For each $x_0 \in X$ and each $\alpha \in \partial X$, 
there exists a unique element $\xi \in \partial_{x_0}X$ 
with $\xi(\infty)=\alpha$.
Thus we may identify $\partial X$ with $\partial_{x_0}X$ for each $x_0 \in X$.

Let $X$ be a proper CAT(0) space 
and $G$ a group which acts on $X$ by isometries.
For each element $g \in G$ and 
each geodesic ray $\xi:[0,\infty)\rightarrow X$, 
a map $g \xi:[0,\infty)\rightarrow X$ 
defined by $(g\xi)(t):=g(\xi(t))$ is also a geodesic ray.
If geodesic rays $\xi$ and $\xi'$ are asymptotic, 
then $g\xi$ and $g\xi'$ are also asymptotic.
Thus $g$ induces a homeomorphism of $\partial X$ and 
$G$ acts on $\partial X$.

A {\it geometric} action on a CAT(0) space 
is an action by isometries which is proper (\cite[p.131]{BH}) 
and cocompact.
We note that every CAT(0) space 
on which a group acts 
geometrically is a proper space (\cite[p.132]{BH}).

Details of CAT(0) spaces and their boundaries are found in 
\cite{BH} and \cite{GH}.

Suppose that a group $G$ acts geometrically on a CAT(0) space $X$.
For a subset $A\subset G$, 
the limit set $L(A)$ of $A$ is defined as 
$L(A)=\overline{Ax_0}\cap \partial X$, 
where $x_0\in X$ and 
$\overline{Ax_0}$ is the closure 
of the orbit $Ax_0$ in $X\cup \partial X$.
We note that 
the limit set $L(A)$ is determined by $A$ and not depend on the point $x_0$.
For $g\in G$, 
we define $F_g$ and $\mathcal{F}_g$ 
as the fixed-point sets of $g$ in $X$ and $\partial X$ respectively, 
that is, 
$F_g=\{x\in X|\, gx=x\}$ and 
$\mathcal{F}_g=\{\alpha\in\partial X|\, g\alpha=\alpha\}$.
By \cite[Corollary~II.2.8~(1)]{BH}, 
we see that $F_g\neq\emptyset$ if and only if 
the order $o(g)$ of $g$ is finite, 
since the action of $G$ on $X$ is proper.
In this paper, we investigate $\mathcal{F}_g$.
We prove the following theorem in Section~2.

\begin{Theorem}\label{Thm1}
Suppose that a group $G$ acts geometrically on a CAT(0) space $X$.
For $g\in G$,
$\mathcal{F}_g=L(Z_g)$, 
where $Z_g$ is the centralizer of $g$, i.e., 
$Z_g=\{v\in G|\, gv=vg\}$.
\end{Theorem}

For a subset $A\subset G$, 
$L(A)\neq\emptyset$ if and only if $A$ is infinite.
Hence Theorem~\ref{Thm1} implies the following corollary.

\begin{Corollary}\label{Cor1}
Suppose that a group $G$ acts geometrically on a CAT(0) space $X$.
For $g\in G$, 
$\mathcal{F}_g\neq\emptyset$ if and only if $Z_g$ is infinite.
\end{Corollary}

We can obtain the following corollary from Corollary~\ref{Cor1}.

\begin{Corollary}\label{Cor2}
Suppose that a group $G$ acts geometrically on CAT(0) spaces $X$ and $X'$.
For $g\in G$, 
$\mathcal{F}_g\neq\emptyset$ if and only if $\mathcal{F}'_g\neq\emptyset$, 
where $\mathcal{F}_g$ and $\mathcal{F}'_g$ are 
the fixed-point sets of $g$ in $\partial X$ and $\partial X'$ respectively.
\end{Corollary}

Let $X$ be a CAT(0) space and let $g$ be an isometry of $X$.
The {\it translation length} of $g$ is the number 
$|g|:=\inf\{d(x,gx)\,|\,x\in X\}$, and 
the {\it minimal set} of $g$ is defined as 
$\Min(g):=\{x\in X\,|\,d(x,gx)=|g|\}$.
An isometry $g$ of $X$ is said to be {\it hyperbolic}, 
if $\Min(g)\neq\emptyset$ and $|g|>0$ (cf.\ \cite{BH}).

We also prove the following theorem in Section~2.

\begin{Theorem}\label{Thm2}
Let $g$ be a hyperbolic isometry of a proper CAT(0) space $X$.
Then the fixed-point set of $g$ in $\partial X$ is 
$\mathcal{F}_g=\partial\Min(g)$, 
where $\partial\Min(g)$ is 
the boundary of the minimal set $\Min(g)$ of $g$.
\end{Theorem}

For a hyperbolic isometry $g$ of a CAT(0) space $X$, 
$\Min(g)$ is the union of the axes of $g$ and 
$\Min(g)$ splits as a product $Y_1\times \Image\sigma$,
where $\sigma:\R\rightarrow X$ is an axis of $g$ (\cite[p.231]{BH}).
Hence 
$\partial\Min(g)$ is the suspension $\partial Y_1*\{g^{-\infty},g^{\infty}\}$.
Theorem~\ref{Thm2} implies that 
the fixed-point set $\mathcal{F}_g$ has a suspension form.
Let $\mathcal{P}_g$ be the periodic-point set of $g$ in $\partial X$.
Then 
$$\mathcal{P}_g=\bigcup_{n\in\N}\mathcal{F}_{g^n}
=\bigcup_{n\in\N}\partial\Min(g^n).$$
Hence the periodic-point set $\mathcal{P}_g$ of $g$ 
is also a suspension form.

We obtain the following corollary from Theorems~\ref{Thm1} and \ref{Thm2}.

\begin{Corollary}\label{Cor3}
Suppose that a group $G$ acts geometrically on a CAT(0) space $X$.
For $g\in G$ such that $o(g)=\infty$,
$\mathcal{F}_g=L(Z_g)=\partial\Min(g)$.
\end{Corollary}

%%%%%%%%%%%%%
% Section 2 %
%%%%%%%%%%%%%
\section{Proof of the main theorems}

We first prove Theorem~\ref{Thm1}.

\begin{proof}[Proof of Theorem~\ref{Thm1}]
Suppose that a group $G$ acts geometrically on a CAT(0) space $X$.
Let $g\in G$ and let $x_0\in X$.

We first show that 
$\mathcal{F}_g\subset L(Z_g)$.
Let $\alpha\in \mathcal{F}_g$ and 
let $\xi:[0,\infty)\rightarrow X$ be the geodesic ray 
such that $\xi(0)=x_0$ and $\xi(\infty)=\alpha$.
Since $\alpha\in \mathcal{F}_g$, $g\alpha=\alpha$.
Hence the geodesic rays $\xi$ and $g\xi$ are asymptotic, 
and there exists a number $M>0$ 
such that $d(\xi(t),g\xi(t))\le M$ for any $t\ge 0$.
Since the action of $G$ on $X$ is cocompact and $X$ is proper, 
$GB(x_0,N)=X$ for some $N>0$.
For each $i\in\N$, 
there exists $v_i\in G$ such that $d(\xi(i),v_ix_0)\le N$.
Then 
\begin{align*}
d(x_0,v_i^{-1}gv_ix_0)&=d(v_ix_0,gv_ix_0) \\
&\le d(v_ix_0,\xi(i))+d(\xi(i),g\xi(i))+d(g\xi(i),gv_ix_0) \\
&\le 2N+M,
\end{align*}
because $d(v_ix_0,\xi(i))\le N$ and $d(\xi(i),g\xi(i))\le M$.
Since the action of $G$ on $X$ is proper, 
the set $\{h\in G|\, d(x_0,hx_0)\le 2N+M\}$ is finite.
Hence there exists $g'\in G$ 
such that $\{i\in \N|\,v_i^{-1}gv_i=g'\}$ is infinite.
Let $\{i_j|\,j\in \N\}=\{i\in \N|\,v_i^{-1}gv_i=g'\}$.
Then $v_{i_j}^{-1}gv_{i_j}=g'$ for each $j\in\N$, and 
$$ v_{i_j}^{-1}gv_{i_j}=g'=v_{i_1}^{-1}gv_{i_1}. $$
Hence $g=(v_{i_j}v_{i_1}^{-1})g(v_{i_1}v_{i_j}^{-1})$.
Thus $v_{i_j}v_{i_1}\in Z_g$ for each $j\in\N$.
Here we note that 
$v_{i_j}v_{i_1}^{-1}\neq v_{i_{j'}}v_{i_1}^{-1}$
if $j\neq j'$.
Since $\{v_{i_j}x_0|\,j\in\N\}$ converges to $\alpha$ in $X\cup\partial X$, 
$\{v_{i_j}v_{i_1}^{-1}x_0|\,j\in\N\}$ converges to $\alpha$.
Hence $\alpha\in L(\{v_{i_j}v_{i_1}^{-1}|\,j\in\N\})\subset L(Z_g)$.
Thus $\mathcal{F}_g\subset L(Z_g)$.

Next, we show that 
$\mathcal{F}_g\supset L(Z_g)$.
Let $\alpha\in L(Z_g)$.
There exists a sequence $\{v_i|\,i\in\N\}\subset Z_g$ 
such that $\{v_ix_0|\,i\in\N\}$ converges to $\alpha$ in $X\cup\partial X$.
For each $i\in\N$, $gv_i=v_ig$.
Here the sequence $\{gv_ix_0|\,i\in\N\}$ converges to $g\alpha$ 
and $\{v_igx_0|\,i\in\N\}$ converges to $\alpha$.
Thus $g\alpha=\alpha$, i.e., 
$\alpha\in\mathcal{F}_g$.

Therefore $\mathcal{F}_g=L(Z_g)$.
\end{proof}

Next we prove Theorem~\ref{Thm2}.

\begin{proof}[Proof of Theorem~\ref{Thm2}]
Let $g$ be a hyperbolic isometry of a proper CAT(0) space $X$ 
and let $x_0\in \Min(g)$.

We first show that $\mathcal{F}_g\subset\partial\Min(g)$.
Let $\alpha\in \mathcal{F}_g$ and 
let $\xi:[0,\infty)\rightarrow X$ be the geodesic ray 
such that $\xi(0)=x_0$ and $\xi(\infty)=\alpha$.
Since $\alpha\in \mathcal{F}_g$, $g\alpha=\alpha$.
Hence the geodesic rays $\xi$ and $g\xi$ are asymptotic.
Now $d(x_0,gx_0)=|g|$ because $x_0\in \Min(g)$.
Hence 
$$d(\xi(t),g\xi(t))\le d(\xi(0),g\xi(0))=d(x_0,gx_0)=|g|.$$ 
This means that $d(\xi(t),g\xi(t))=|g|$ and $\xi(t)\in \Min(g)$ 
for each $t\ge 0$.
Hence $\Image\xi\subset\Min(g)$ and $\alpha\in\partial\Min(g)$.
Thus $\mathcal{F}_g\subset\partial\Min(g)$.

Next we show that $\mathcal{F}_g\supset\partial\Min(g)$.
Let $\alpha\in \partial\Min(g)$ and 
let $\xi:[0,\infty)\rightarrow X$ be the geodesic ray 
such that $\xi(0)=x_0$ and $\xi(\infty)=\alpha$.
Then $\Image\xi\subset\Min(g)$, 
since $\alpha\in\partial\Min(g)$ and $\Min(g)$ is convex.
Hence $\xi(t)\in\Min(g)$ and $d(\xi(t),g\xi(t))=|g|$ for any $t\ge0$.
This means that $\xi$ and $g\xi$ are asymptotic and $\alpha=g\alpha$, 
i.e., $\alpha\in\mathcal{F}_g$.
Thus $\mathcal{F}_g\supset\partial\Min(g)$.

Therefore $\mathcal{F}_g=\partial\Min(g)$.
\end{proof}

%%%%%%%%%%%%%
% Section 4 %
%%%%%%%%%%%%%
\section{Remarks}

Let $g$ be a hyperbolic isometry of a proper CAT(0) space $X$.
For each $n\in \N$, 
$\mathcal{F}_{g^n}=\partial\Min(g^n)$ by Theorem~\ref{Thm2}.
We note that an axis of $g$ is also an axis of $g^n$ for each $n\in\N$.
Then the minimal set $\Min(g^n)$ of $g^n$ 
splits as a product $\Min(g^n)=Y_n\times \Image\sigma$, 
where $\sigma:\R\rightarrow X$ is an axis of $g$ (cf.\ \cite{BH}).
Here $Y_n\subset Y_{kn}$ for each $n,k\in\N$.
The fixed-point set of $g$ in $\partial X$ is 
$$ \mathcal{F}_{g}=\partial\Min(g)=\partial(Y_1\times\Image\sigma)
=\partial Y_1*\{g^{\infty},g^{-\infty}\}.$$
Also the periodic-point set of $g$ is 
\begin{align*}
\mathcal{P}_g&=\bigcup_{n\in\N}\mathcal{F}_{g^n}=\bigcup_{n\in\N}\partial\Min(g^n)\\
&=\bigcup_{n\in\N}\partial(Y_n\times\Image\sigma)
=(\bigcup_{n\in\N}\partial Y_n)*\{g^{\infty},g^{-\infty}\}.
\end{align*}
Thus the fixed-point set and the periodic-point set of $g$ have 
suspension forms.

Let $A$ be the union of geodesic lines 
which are parallel to an axis $\sigma$ of $g$.
Then $A$ splits as a product $Y\times\Image\sigma$.
By the above argument, 
$\mathcal{P}_g=\bigcup_{n\in\N}\partial\Min(g^n)\subset \partial A$ and 
$\bigcup_{n\in\N}\partial Y_n\subset \partial Y$.

Here the following problem arises.

\begin{Problem}
Is it always the case that 
$\mathcal{P}_g=\partial A$ 
(i.e.\ $\bigcup_{n\in\N}\partial Y_n=\partial Y$)?
\end{Problem}

Let $g$ be a hyperbolic isometry of a Gromov hyperbolic space $X$.
By an easy argument, 
we see that 
the fixed-point set $\mathcal{F}_g$ 
of $g$ in the boundary $\partial X$ is 
the two-points set $\{g^{\infty},g^{-\infty}\}$.
Also for each $\alpha\in\partial X\setminus\{g^{-\infty}\}$, 
the sequence $\{g^i\alpha\,|\,i\in\N\}$ converges to $g^\infty$ in $\partial X$.

Here the following problem arises.

\begin{Problem}
Let $g$ be a hyperbolic isometry of a proper CAT(0) space $X$ and 
let $A$ be the union of geodesic lines 
which are parallel to an axis $\sigma$ of $g$.
Is it the case that 
for each $\alpha\in\partial X\setminus\partial A$, 
the sequence $\{g^i\alpha\,|\,i\in\N\}$ converges to $g^\infty$ 
in $\partial X$?
\end{Problem}

Also the following problem arises.

\begin{Problem}
Suppose that a group $G$ acts geometrically on a CAT(0) space $X$.
Is it always the case that 
there do not exist $g\in G$ and $\alpha,\beta\in\partial X$ 
such that 
${\limsup}_{n\rightarrow\infty}d_{\partial X}(g^n\alpha,g^n\beta)>0$ 
and 
${\liminf}_{n\rightarrow\infty}d_{\partial X}(g^n\alpha,g^n\beta)=0$?
\end{Problem}

%%%%%%%%%%%%%%%%%%%%%%%%%%%%%%%%%%%%%
%             REFERENCES            %
%%%%%%%%%%%%%%%%%%%%%%%%%%%%%%%%%%%%%
%

%
\end{document}